\documentclass{article}
\usepackage{amsfonts}
\usepackage{amssymb}
\title{\large \bf Automorphisms of groups and converse of Schur's theorem}
\author{\small \bf Deepak Gumber, Hemant Kalra and Sandeep Singh \\
\small \em School of Mathematics and Computer Applications\\
\small \em Thapar University, Patiala - 147 004,
India}

\date{}
\newtheorem{thm}{Theorem}[section]
\newtheorem{lm}[thm]{Lemma}

\newtheorem{cor}[thm]{Corollary}

\begin{document}

\maketitle
\begin{abstract}
An automorphism of a group $G$ is called an IA-automorphism if it induces the identity automorphism on the abelianized group $G/G^{\prime}$. Let $\mathrm{IA}(G)$ denote the group of all IA-automorphisms of $G$.
We classify all finitely generated nilpotent groups $G$ of class 2 for which $\mathrm{IA}(G)\simeq \mathrm{Inn}(G)$. In particular, we classify all finite nilpotent groups of class 2 for which each $\mathrm{IA}$-automorphism is inner. As consequences, we give  surprisingly very easy proofs of converse of Schur's theorem and also prove some other related results.
\end{abstract}

\vspace{2ex}

\noindent {\bf 2010 Mathematics Subject Classification:}
20F28, 20F18.

\vspace{2ex}

\noindent {\bf Keywords:} IA-automorphism, Schur's theorem, nilpotent group.

\section{\normalsize Introduction}

Let $G$ be any group and let $G^{\prime}$ and $Z(G)$ respectively denote the derived group and the center of $G$. By $\mathrm{Inn}(G)$ we denote the inner automorphism group of $G$. Following Bachmuth [2], we call an automorphism of $G$ an IA-automorphism if it induces the identity automorphism on the abelianized group $G/G^{\prime}$. Let $\mathrm{IA}(G)$ denote the group of all IA-automorphisms of $G$ and let $\mathrm{IA}(G)^*$ denote its subgroup consisting of those IA-automorphisms which fix $Z(G)$ element-wise. IA-automorphisms have been well studied in past (see for example [2], [4], [5], [9], [11]). If $G$ is a free group of rank 2, then it is a classical result of Nielsen [9] that $\mathrm{IA}(G)= \mathrm{Inn}(G)$. Bachmuth [2] obtained the same result for free metabelian groups of rank 2. Recently Attar [1, Theorem 2.1] has proved that if $G$ is a finite $p$-group of class 2, then $\mathrm{IA}(G)= \mathrm{Inn}(G)$ if and only if $G^{\prime}$ is cyclic and $\mathrm{IA}(G)=\mathrm{IA}(G)^*$. In section 2, we classify all finitely generated nilpotent groups $G$ of class 2 for which (i) $\mathrm{IA}(G)\simeq \mathrm{Inn}(G)$ and (ii) $\mathrm{IA}(G)^*\simeq \mathrm{Inn}(G)$. In particular,  we classify all finite nilpotent groups $G$ of class 2 for which (i) $\mathrm{IA}(G)= \mathrm{Inn}(G)$ and (ii) $\mathrm{IA}(G)^*=\mathrm{Inn}(G)$.

In the year 1904 Schur [12] proved the result: If the index of the center of a group is finite, then its derived group is also finite. We call this result as Schur's theorem. Neumann [8, Corrolary 5.41] proved the converse of Schur's theorem for finitely generated groups (the converse of Schur's theorem is not true in general for example for an infinite extra-special $p$-group for an odd prime $p$). More precisely, he proved that if $G$ is finitely generated and if its
derived group $G^{\prime}$ has finite order, then the index $[G:Z(G)]$ of the centre  is also finite. We call this result as Neumann's theorem. Niroomand [10]  generalized the Neumann's theorem to the following: If $G^{\prime}$ is finite and $G/Z(G)$ is finitely generated by $d$ elements, then 
\begin{equation}
|G/Z(G)|\le |G^{\prime}|^d.
\end{equation}
Sury [13] generalized Niroomand's result to the following: Let $G$ be a group such that the set $K(G)$ of its commutators is finite. Then $G^{\prime}$ is finite and  if $G/Z(G)$ can be generated by $d$ elements, then
\begin{equation}
|G/Z(G)| \le |K(G)|^d.
\end{equation}
Recently Yadav [15, Theorem A] has further generalized Niroomand's and Sury's results and has raised the question: what can one say about non-abelian groups for which equality holds in (1)? He himself classifies, upto isoclinism, all nilpotent groups $G$ of class 2 such that $G/Z(G)$ is finite and equality holds in (1). It can be easily seen that if $G$ is a finite $p$-group of maximal class, then equality holds in (1) if and only if $|G|=p^3$. In section 3, we prove that if $G$ is a finite $p$-group of co-class 2 for which equality holds in (1), then $|G|=p^4$ or $p^5$; and give necessary and sufficient conditions for a finitely generated nilpotent group $G$ of class 2 for which $G/Z(G)\simeq {G^{\prime}}^{r(G/Z(G))}$. As a consequence of the fact that $\mathrm{Inn}(G)\le\mathrm{IA}(G)^*$, we also obtain the converses of Schur's theorem proved by Niroomand [10], Sury [13] and Yadav [15, Theorem A]. Quite unexpectedly, our proofs are easy, in-fact very easy.  Yadav [15, Question 2] has also raised the question that if $G$ is a finite nilpotent group for which the equality holds in (1), then is it necessary that $G^{\prime}$ is cyclic? We answer the question in negative and give a counter example of minimal possible order.

By $\mathrm{Hom}(G,A)$ we denote the group of all homomorphisms of
$G$ into an abelian group $A$. For $x\in G$, the set $\{[x,g]|g\in G\}$ of commutators is denoted as $[x,G]$. The frattini subgroup of $G$ is denoted as $\Phi(G)$ and $(G)_{ab}$ denotes the abelianized group $G/G^{\prime}$.  By $C_{p}$ we denote a cyclic group of order $p$ and by $X^n$ we denote the direct product of $n$-copies of a group $X$. The rank of a group $G$ is the smallest cardinality of a generating set of $G$. The rank, torsion rank and torsion-free rank of a group $G$ are respectively denoted as $r(G)$, $d(G)$ and $\rho(G)$. For a finite group $G$, $G_p$, $\mathrm{exp}(G)$ and $\pi(G)$ respectively denote the Sylow $p$-subgroup, exponent and the set of primes dividing the order of $G$. All other unexplained notations, if any, are standard. We shall use the  following two well known lemmas very frequently without any further referring.

\begin{lm}
Let $U,V$ and $W$ be abelian groups. Then\\
$(i)$ if $U$ is torsion-free of rank $m$, then $\mathrm{Hom}(U, V)\simeq V^{m}$,\\ 
$(ii)$ if $U$ is torsion and $V$ is torsion-free, then 
$\mathrm{Hom}(U, V)=1$, and\\
$(iii)$ if $U$ and $V$ are finite and $\mathrm{exp}(U)=\mathrm{exp}(V)$, then $Hom(U, V)\simeq U$ if and only if $V$ is a cyclic group. 
\end{lm}

\begin{lm} Let $G$ be a finitely generated nilpotent group of class $2$. Then $\mathrm{exp}(T(G/Z(G)))=\mathrm{exp}(T(G^{\prime}))$, where $T(H)$ denotes the torsion part of a group $H$.
\end{lm}

Let $G$ be a group and $X$ be a central subgroup of $G$. If $\theta\in\mathrm{Hom}((G/X)_{ab},X)$, then $T_{\theta}:G\rightarrow G$ defined as $T_{\theta}(g)=g\theta(\overline{g})$ is an automorphism of $G$, where $\overline{g}$ is the coset determined by $g$ in $(G/X)_{ab}$. It is clear that $\theta \mapsto T_{\theta}$ is
an isomorphism between $\mathrm{Hom}((G/X)_{ab},X)$ and the group of all automorphisms of $G$ that induce the identity on both $X$ and $G/X$. This observation is from [4, Sec. 2] and is a little specific one. We can, in-fact, have the following more general result:

\begin{lm}
Let $G$ be any group and $Y$ be a central subgroup of $G$ contained in a normal subgroup $X$ of $G$. Then the group of all automorphisms of $G$ that induce the identity on both $X$ and $G/Y$ is isomorphic onto
$\mathrm{Hom}(G/X,Y)$.
\end{lm}

The following is an easy particular consequence:

\begin{lm}
Let $G$ be a nilpotent group of class $2$. Then $$\mathrm{IA}(G)\simeq \mathrm{Hom}(G/G^{\prime},G^{\prime})\;\;\mathrm{and} \;\;\mathrm{IA}(G)^*\simeq \mathrm{Hom}(G/Z(G),G^{\prime}).$$ 
\end{lm}

\section{\normalsize IA-automorphisms}

Let $G$ be a finitely generated nilpotent group of class 2. Let $G/Z(G)\simeq A\times \mathbb{Z}^a$, $G/G^{\prime}\simeq B\times \mathbb{Z}^b$ and $G^{\prime}\simeq C\times \mathbb{Z}^c$, where $A, B,C$  are respective torsion parts and $a,b,c$ are respective torsion-free ranks of $G/Z(G),G/G^{\prime}$ and $G^{\prime}$. Let $\pi(A)=\{p_1,p_2,\ldots,p_d\}$ and $\pi(B)=\{p_1,p_2,\ldots,p_e\}$, where $d\le e$ because $G/Z(G)$ is a quotient group of $G/G^{\prime}$.
Suppose that
 $$A_{p_{i}}\simeq \displaystyle\prod_{j=1}^{m_{i}}C_{p_{i}^{\alpha_{ij}}}\;\; 
\mathrm{and}\;\;B_{p_{i}}\simeq \displaystyle\prod_{j=1}^{n_{i}}C_{p_{i}^{\beta_{ij}}}.$$
Then
\[
A\simeq  \displaystyle\prod_{i=1}^d
\displaystyle\prod_{j=1}^{m_{i}}C_{p_{i}^{\alpha_{ij}}}\;\;
\mathrm{and}\;\;B\simeq 
\displaystyle\prod_{i=1}^e
\displaystyle\prod_{j=1}^{n_{i}}C_{p_{i}^{\beta_{ij}}},
\]
where $ \alpha_{11}\geq \alpha_{12}\geq\ldots\ge\alpha_{1m_{1}},\alpha_{21}\geq \alpha_{22}\geq\ldots\geq \alpha_{2m_{2}},\ldots,\alpha_{d1}\geq \alpha_{d2}\geq\ldots\geq\alpha_{dm_{d}}$ and $ \beta_{11}\geq \beta_{12}\ge\ldots\ge\beta_{1n_{1}}, \beta_{21}\geq \beta_{22}\geq\ldots\beta_{2n_{2}},\ldots,\beta_{e1}\geq \beta_{e2}\geq\ldots\geq \beta_{en_{e}}$ are positive integers. Since $G/Z(G)$ is a quotient group of $G/G^{\prime}$, it follows that $m_i\le n_i$ and $\alpha_{ij}\le\beta_{ij}$ for all $i,1\le i\le d$ and for all $j,1\le j\le m_i$. We fix these notations, unless or otherwise stated, for the rest of this paper.

\begin{thm}
Let $G$ be a finitely generated nilpotent group of class $2$. Then \\
$(i)$ if $G$ is torsion-free, then $\mathrm{IA}(G)\simeq \mathrm{Inn}(G)$ if and only if  $G^{\prime}$ is cyclic and $\rho (G/Z(G))=\rho(G/G^{\prime})$.\\
$(ii)$ if $G$ is not torsion-free, then $\mathrm{IA}(G)\simeq \mathrm{Inn}(G)$ if and only if one of the following four conditions holds:\\
$(a)$ $G$ is finite, $m_{i}=n_{i}$ for $1\leq i \leq d$, $G^{\prime}$ is cyclic, and for each $i,\;1\leq i \leq d$, either $(G/Z(G))_{p_{i}} $ is homocyclic or $\alpha_{i1}=\alpha_{it_{i}}$ for $ i1\leq it_{i}\leq i(r_{i}-1)$ and $\beta_{it_{i}}= \alpha_{it_{i}}$ for $ir_{i}\leq it_{i}\leq in_{i}$, where $ir_{i}$ is the smallest positive integer between $i1$ and $ in_{i}$ such that $\beta_{ir_{i}}< \alpha_{i1}$.\\
$(b)$ $G^{\prime}\simeq \mathbb{Z}$ and $\rho(G/Z(G))=\rho(G/G^{\prime})$.\\ 
$(c)$ $G^{\prime}$ is torsion, $G/G^{\prime}$ is torsion-free and $A\simeq C^b.$\\
$(d)$ $G/Z(G)$ and $G^{\prime}\simeq C\times \mathbb{Z}$ are mixed groups, $G/G^{\prime}$ is torsion-free, $A\simeq C^b$ and $\rho (G/Z(G))= \rho (G/G^{\prime})$.\\
$(iii)$ $\mathrm{IA}(G)^*\simeq \mathrm{Inn}(G)$ if and only if $G^{\prime}$ is cyclic.
\end{thm}
\noindent{\bf Proof.} (i) First Suppose that $\mathrm{IA}(G)\simeq \mathrm{Inn}(G)$. Then
$$\mathrm{IA}(G)\simeq \mathrm{Hom}(G/G^{\prime},G^{\prime})\simeq G/Z(G)$$ 
by Lemma 1.4. If $G/G^{\prime}$ is finite, then $G/Z(G)$ is finite and hence by Schur's theorem $G^{\prime}$ is finite, which is not so. Thus $G/G^{\prime}\simeq B\times \mathbb{Z}^b$, where $b>0$. If $G/Z(G)$ has a torsion element, then so does $G^{\prime}$, which is not so. Thus $G/Z(G)$ is torsion-free. It then follows that $\mathbb{Z}^{bc}\simeq  \mathbb{Z}^a$.  Since $a\le b$, $c=1$ and $a=b$. Hence $G^{\prime}$ is cyclic and $\rho(G/Z(G))=\rho(G/G^{\prime})$. The converse follows easily.\\
(ii) It is not very hard to see that if $G$ satisfies any of the four conditions, then $\mathrm{IA}(G)\simeq \mathrm{Inn}(G)$. Conversely, if $\mathrm{IA}(G)\simeq \mathrm{Inn}(G)$, then
\begin{equation}
\mathrm{Hom}(B\times \mathbb{Z}^b,C\times \mathbb{Z}^c)\simeq A\times \mathbb{Z}^a.
\end{equation}
First suppose that $G$ is finite. Then $\mathrm{Hom}(B,C)\simeq A$. Since $d(A)\le d(B)$ and $\mathrm{exp}(A)=\mathrm{exp}(C)$, $G^{\prime}$ is cyclic and thus  $G^{\prime}\simeq C_{{p_{1}}^{\alpha_{11}}}\times C_{{p_{2}}^{\alpha_{21}}}\times \ldots \times C_{{p_{d}}^{\alpha_{d1}}}$. Observe that
\[\begin{array}{lcl}
\mathrm{Hom}(B,C)&\simeq&\mathrm{Hom}(\displaystyle\prod_{i=1}^e\displaystyle\prod_{j=1}^{n_{i}}
C_{p_{i}^{\beta_{ij}}}, \displaystyle\prod_{i=1}^{d}C_{p_{i}^{\alpha_{i1}}})\\
&\simeq&\mathrm{Hom}(\displaystyle\prod_{i=1}^d\displaystyle\prod_{j=1}^{n_{i}}
C_{p_{i}^{\beta_{ij}}}, \displaystyle\prod_{i=1}^{d}C_{p_{i}^{\alpha_{i1}}})\\
&\simeq&\displaystyle\prod_{i=1}^d\mathrm{Hom}(\displaystyle\prod_{j=1}^{n_{i}}
C_{p_{i}^{\beta_{ij}}}, \displaystyle\prod_{i=1}^{d}C_{p_{i}^{\alpha_{i1}}})\\
&\simeq&\displaystyle\prod_{i=1}^d\mathrm{Hom}(\displaystyle\prod_{j=1}^{n_{i}}
C_{p_{i}^{\beta_{ij}}}, C_{p_{i}^{\alpha_{i1}}}),\\
\end{array}\]
and $A\simeq\displaystyle\prod_{i=1}^d\displaystyle\prod_{j=1}^{m_{i}}
C_{p_{i}^{\alpha_{ij}}}$. It thus follows that for each $i,1\le i\le d$,
$$\mathrm{Hom}(\displaystyle\prod_{j=1}^{n_{i}}
C_{p_{i}^{\beta_{ij}}}, C_{p_{i}^{\alpha_{i1}}})\simeq \displaystyle\prod_{j=1}^{m_{i}}
C_{p_{i}^{\alpha_{ij}}}.$$
This immediately implies that $m_i=n_i$. Now either $\beta_{ij}\ge \alpha_{i1}$ for each $j,1\le j\le n_i$, or there exist smallest positive integer $ir_{i}$ between $i1$ and $in_i$ such that $\beta_{ir_{i}}< \alpha_{i1}$. In the first case $(G/Z(G))_{p_{i}}$ is homocyclic, and in the second case
$\alpha_{i1}=\alpha_{it_{i}}$ for $i1\le it_{i}\le i(r_{i}-1)$ and $\beta_{it_{i}}=\alpha_{it_{i}}$ for $ir_i \le it_{i}\le in_i$.

Next suppose that $G$ is infinite. The derived group $G^{\prime}$ can be torsion, mixed, or torsion-free. First assume that $G^{\prime}$ is torsion-free. Then $G/Z(G)$ is torsion-free and $G/G^{\prime}$ is not torsion by Schur's theorem. Equation (3) then implies that $\mathbb{Z}^{bc}\simeq \mathbb{Z}^a$. It follows that $a=b$ and $c=1$. Next assume that $G^{\prime}$ is torsion. Then $G/Z(G)$ is torsion by Neumann's theorem. If $G/G^{\prime}$ is torsion, then $G$ is finite. If $G/G^{\prime}$ is torsion free, then (3) implies that $A\simeq C^b$; and if $G/G^{\prime}$ is a mixed group, then (3) implies that $A\simeq \mathrm{Hom}(B,C)\times C^b$. Since $d(A)\le d(B)$, $b=0$ and hence $G$ is finite, which is not so.   Finally assume that $G^{\prime}$ is a mixed group. Then $G/G^{\prime}$ is either mixed or torsion-free by Schur's theorem. 
If $G/G^{\prime}$ is torsion-free, then (3) implies that $A\simeq C^b$, $a=b$, $c=1$; and if $G/G^{\prime}$ is mixed, then (3) implies that $A\simeq \mathrm {Hom}(B,C)\times C^b$ and $a=bc$. It follows that $a=b=0$. Therefore both $G/Z(G)$ and $G/G^{\prime}$ are finite. By Schur's theorem $G^{\prime}$ is finite and hence $G$ is finite, a contradiction. \\
(iii) If $G^{\prime}$ is cyclic, finite or infinite, then, by using Lemma 1.4, it is easy to see that $\mathrm{IA}(G)^*\simeq \mathrm{Inn}(G)$. To prove the converse, first suppose that $G$ is torsion-free. Then $G^{\prime}$ is torsion-free and hence $G/Z(G)$ is torsion-free. It follows that  $\mathrm{Hom}(\mathbb{Z}^a,\mathbb{Z}^c)\simeq \mathbb{Z}^a$, and hence $G^{\prime}\simeq \mathbb{Z}$ is cyclic. 
          
Next suppose that $G$ is torsion. Then both $G/Z(G)$ and $G^{\prime}$ are torsion. Thus $\mathrm{Hom}(A,C)\simeq A$ and hence $G^{\prime}$ is finite cyclic.

Finally suppose that $G$ is a mixed group. Then either $G^{\prime}$ is  a torsion, or torsion-free, or a mixed group. The case that $G^{\prime}$ is a torsion-free group can be handled as before.  If $G^{\prime}$ is torsion, then $G/Z(G)$ is torsion by Neumann's theorem and hence as above $G^{\prime}$ is finite cyclic. Now suppose that $G^{\prime}$ is a mixed group. Then $G/Z(G)$ is either a mixed or a torsion-free group by Schur's theorem. If $G/Z(G)$ is mixed, then  $\mathrm{Hom}(A,C)\times C^a\simeq A$. Since 
$\mathrm{exp}(A)=\mathrm{exp}(C)$, $a=0$ and $G^{\prime}$ is a finite cyclic group. Finally if $G/Z(G)$ is torsion-free, then $\mathrm{Hom}(\mathbb{Z}^a, C\times \mathbb{Z}^c)\simeq \mathbb{Z}^a$ implies that $G^{\prime}\simeq \mathbb{Z}$.
\hfill $\Box$\\

Let $\mathrm{Aut}_c(G)$ denote the group of all conjugacy class preserving automorphisms of a group $G$. 
As a consequence of Theorem 2.1(iii), which generalizes corollary 2.3 of Attar [1], we have the following small but significant generalization of corollary 3.6 of Yadav [14].

\begin{cor}
Let $G$ be a finitely generated nilpotent group of class $2$ such that $G^{\prime}$ is cyclic and $G/Z(G)$ is finite. Then $\mathrm{Aut}_c(G)=\mathrm{Inn}(G)$. 
\end{cor}
{\bf Proof.} The proof follows from the fact that
$\mathrm{Inn}(G)\le\mathrm{Aut}_c(G)\le\mathrm{IA}(G)^*.$ \hfill $\Box$\\

\begin{cor}
Let $G$ be a $2$-generated finite nilpotent group of class $2$. Then any $\mathrm{IA}$-automorphism of $G$ is an inner automorphism. 
\end{cor}
{\bf Proof.} Observe that $G^{\prime}$ is cyclic and $d(G/Z(G))=d(G/G^{\prime})=2$. Also $G/Z(G)$ is homocyclic by [7, Lemma 0.4]. The result now follows from  Theorem 2.1. \hfill $\Box$

\begin{cor}
Let $G$ be a $2$-generated torsion-free nilpotent group of class $2$. Then $\mathrm{IA}(G)\simeq \mathrm{Inn}(G)$.
\end{cor}
{\bf Proof.} Observe that $G^{\prime}\simeq \mathbb{Z}$ and $\rho(G/Z(G))=\rho(G/G^{\prime})=2$. Thus  $\mathrm{IA}(G)\simeq \mathrm{Inn}(G)$  by Theorem 2.1.\hfill $\Box$ \\

Let $G$ be a finite $p$-group of nilpotency class 2 and let
$G/Z(G)\simeq\prod_{i=1}^t C_{p^{\alpha_i}}$ and
$G/G^{\prime}\simeq\prod_{j=1}^s C_{p^{\delta_j}}$,
where  $\alpha_1\geq\alpha_2\geq\ldots\geq\alpha_t$ and $\delta_1\geq\delta_2\geq\ldots\geq\delta_s$ are positive integers. Since $G/Z(G)$ is a quotient group of $G/G^{\prime}$, $t\le s$ and $\alpha_i\le \delta_i$ for $1\le i\le t$. Attar [1, Theorem 2.1] has proved that if $G$ is a finite $p$-group of class 2, then $\mathrm{IA}(G)= \mathrm{Inn}(G)$ if and only if $G^{\prime}$ is cyclic and $\mathrm{IA}(G)=\mathrm{IA}(G)^*$. Our next corollary, which is a particular case of Theorem 2.1, gives an explicit interpretation of this result.

\begin{cor}
Let $G$ be a finite $p$-group of class $2$. Then $\mathrm{IA}(G)= \mathrm{Inn}(G)$ if and only if $G^{\prime}$ is cyclic, $d(G/Z(G))=d(G/G^{\prime})$, and either $G/Z(G)$ is homocyclic or $\alpha_i=\alpha_1$ for $1\le i\le r-1$ and $\delta_i=\alpha_i$ for $r\le i\le s$, where $r$ is the smallest positive integer such that $\delta_r< \alpha_1.$ 
\end{cor}

\section {\normalsize Converse of Schur's Theorem}

Let $G$ be an arbitrary group such that $G^{\prime}$ is finite and $G/Z(G)$ is finitely generated. Let $G/Z(G)$ be minimally  generated by $x_1Z(G),x_2Z(G),\ldots,x_dZ(G)$ and let $\alpha\in \mathrm{IA}(G)^*$. Then $\alpha$ fixes $Z(G)$ element-wise and $\alpha(x_i)=x_iy_i$, where $y_i\in G^{\prime}$ for each $i,1\le i\le d$. It thus follows that $\mathrm{IA}(G)^*$ is finite and $|\mathrm{IA}(G)^*|\le |G^{\prime}|^d.$ We have thus proved the following.
\begin{thm}
Let $G$ be an arbitrary group such that $r(G/Z(G))$ and $G^{\prime}$ are
finite. Then $|\mathrm{IA}(G)^*|\le |G^{\prime}|^{r(G/Z(G))}$.
\end{thm}

In particular, since $\mathrm{Inn}(G)\le \mathrm{IA}(G)^*$, we obtain the following main theorem of Niroomand [10].

\begin{cor}
Let $G$ be an arbitrary group such that $r(G/Z(G))$ and $G^{\prime}$ are
finite. Then $|G/Z(G)|\le |G^{\prime}|^{r(G/Z(G))}$.
\end{cor}

Now, rather than $G^{\prime}$, suppose that $[x_i,G]$ is finite for all $i,1\le i\le d$. Let $a\in G$ and let $\alpha$ be the inner automorphism of $G$ defined by conjugation by $a$. Then $\alpha(x_i)=a^{-1}x_ia=x_i[x_i,a]$. Since $|[x_i,G]|=|x_{i}^{G}|\le |K(G)|$, we have thus generalized the main theorem of Sury [13] and proved Theorem A of Yadav [15] in the following theorem.

\begin{thm}
Let $G$ be an arbitrary group such that $G/Z(G)$ is finitely generated by $x_1Z(G),x_2Z(G),\ldots ,x_dZ(G)$ and $[x_i,G]$ is finite for all $i,1\le i\le d$. Then $|G/Z(G)|\le \Pi_{i=1}^{d}[x_i,G]$ and $G^{\prime}$ is finite.
\end{thm}

\begin{thm}
Let $G$ be an arbitrary group. Then $\mathrm{Inn}(G)$ is finite if and only if  $\mathrm{IA}(G)^*$ is finite.
\end{thm}
{\bf Proof.} If $\mathrm{Inn}(G)$ is finite, then $G'$ is finite by Schur's theorem. It then follows from Theorem 3.1 that $\mathrm{IA}(G)^*$ is finite. Converse is obvious.
 \hfill$\Box$\\

Niroomand [10] claims that inequality (1) is sharp and the bound he obtains is best possible. He gives example of the group $Q_8$, the  quaternion group of order 8, to support his claim. It is clear from Theorem 3.3 that bound obtained by Niroomand is not the best possible. The misinterpretation was because of the group $G=Q_8$ in which $G^{\prime}=[x,G]$ for all $x\in G-\Phi(G)$. We observe from inequalities (1), (2) and Theorem 3.3 that 
$$|G/Z(G)|\le \Pi_{i=1}^{d}|[x_i,G]|\le |K(G)|^d\le |G^{\prime}|^d.$$
It follows that if $G$ is an arbitrary group such that $G/Z(G)$ is finitely generated, $G^{\prime}$ is finite and equality holds in (1), then $G^{\prime}=K(G)$. In particular, if $G$ is a finite group for which equality holds in (1), then $G'=K(G)$. This gives a new tool for finding groups $G$ for which $G'=K(G)$ (for more details on this problem one can see the excellent survey article of Kappe and Morse [6]). Also observe that if the equality is achieved in (1) by a group $G$, then for that group $\mathrm{Inn}(G)=\mathrm{Aut}_c(G)=\mathrm{IA}(G)^*$. It now becomes interested to classify the groups for which equality holds in (1). Of course it holds for all abelian groups. Yadav [15] also has asked the question that for what non-abelian groups equality holds in (1). He himself has classified, upto isoclinism, all such nilpotent groups $G$ of class 2 with $G/Z(G)$ finite. We here prove the following.

\begin{thm}
Let $G$ be a finitely generated nilpotent group of class $2$. Then
$$G/Z(G)\simeq G^{\prime^{r(G/Z(G))}}$$
if and only if $G^{\prime}$ is cyclic and $G/Z(G)$ is homocyclic.
\end{thm}
{\bf Proof.} First suppose that $G^{\prime}$ is cyclic and $G/Z(G)$ is homocyclic. If $G^{\prime}$ is torsion-free, then $G/Z(G)$ is also torsion-free. It follows that
$$G/Z(G)\simeq \mathbb{Z}^{r(G/Z(G))}\simeq G^{\prime^{r(G/Z(G))}}.$$
In case $G^{\prime}$ is finite, then $G/Z(G)$ is also finite by Neumann's theorem. Since exponents of $G/Z(G)$ and $G^{\prime}$ are same, $G/Z(G)\simeq G^{\prime^{r(G/Z(G))}}$.

Conversely suppose that $G/Z(G)\simeq G^{\prime^{r(G/Z(G))}}$. It is sufficient to prove that $G^{\prime}$ is cyclic. If $G'$ is finite, then $G'$ is cyclic by Theorem 2.1(iii). If $G^{\prime}$ is torsion-free, then $G/Z(G)$ is torsion-free. Thus $G^{\prime}\simeq \mathbb{Z}^c$ and $G/Z(G)\simeq \mathbb{Z}^a$. Then $\mathbb{Z}^a\simeq \mathbb{Z}^{ac}$ implies that $G^{\prime}\simeq \mathbb{Z}$ is cyclic.
 Finally suppose that $G'=C\times \mathbb{Z}^c$ is mixed. Then by assumption $A\times \mathbb{Z}^a\simeq (C\times\mathbb{Z}^c)^{r(G/Z(G))}$. This is not possible because $a<r(G/Z(G))$ and $cr(G/Z(G))\ge r(G/Z(G))$.
\hfill$\Box$

\begin{thm}
Let $G$ be a finite non-abelian $p$-group of co-class $2$ for which equality holds in $(1)$. Then $|G|=p^4$ or $p^5$.
\end{thm}
{\bf Proof.} Let $|G|=p^n$. Then $p\le |Z(G)|\le p^2$ and $p^{n-3}\le |G^\prime| \le p^{n-2}$. Let $C^*$ denote the group of all automorphisms of $G$ which centralize both $G/Z(G)$ and $Z(G)$. It follows by assumption that $\mathrm{IA}(G)^*=\mathrm{Inn}(G)$ and by Lemma 1.3 that $C^*\approx \mathrm{Hom}(G/G^\prime, Z(G))$. Observe that if $Z(G)\le G^\prime$, then $C^*\le \mathrm{IA}(G)^*=\mathrm{Inn}(G)$ and thus $C^*=Z(\mathrm{Inn}(G))$. Since $d(G)\ge 2$, $Z(\mathrm{Inn}(G))$ cannot be cyclic. First suppose that $|Z(G)|=p$. Then $|Z_2(G)|=p^3$ and $|C^*|=|Z(\mathrm{Inn}(G))|=p^2$. Thus $d(G)=d(G/Z(G))=2$. If $|G^\prime|=p^{n-2}$, then by (1) $n=3$, which is not possible; and if $|G^\prime|=p^{n-3}$, then by (1) $n=5$. Next suppose that $|Z(G)|=p^2$. Then  $Z(G)$ is not contained in $G^\prime$, because then $Z(\mathrm{Inn}(G))$ is cyclic. If $|G^\prime|=p^{n-2}$, then $G'=\Phi(G)$ and thus $d(G/Z(G))=1$, a contradiction. Therefore $G'=p^{n-3}$, $d(G/Z(G))=2$ and hence $n=4$ by (1). \hfill $\Box$

We end up the section by giving a counter example of minimal order to the following question posed by Yadav [15, Question 2]: Let $G$ be a finite nilpotent group for which equality holds in (1). Is it true that $G^{\prime}$ is cyclic? It follows from Theorem 2.1(iii) that if there is a finite $p$-group $G$ for which equality holds in (1) but $G^{\prime}$ is not cyclic, then it cannot be of class 2 or of maximal class and hence its order must be bigger than $p^4$. Our group is of order $32$ with elementary abelian derived group.  \\

\noindent {\bf Example:} Consider the group
$$G=\langle x,y| x^2y^{-4}=[x,y,x]=[x,y,y]y^{-4}=1\rangle$$
which is of nilpotency class 3. Since $d(G)=d(G/Z(G))=2$, $Z(G)\le \Phi(G)$. Take $u=[x,y]$. It is easily seen that every element of $G$ is of the form $x^iy^ju^k$, $1\le i,k \le 2$, $1\le j\le 8$; and $|x|=4, |y|=8$, $|u|=2$.
Since $|\Phi(G)|=8$, $\Phi(G)=\langle y^2, u\rangle$. It is easy to see that $Z(G)=\langle y^4\rangle$ and $|Z(G)|=2$. If $|G'|=8$, then $|\gamma_3(G)|=4$ by [3, Theorem 1.5], a contradiction to $\gamma_3(G)\le Z(G)$. Thus $|G^{\prime}|=4$. Since $y^4,u\in G'$ are of order 2, $G^{\prime}$ is elementary abelian. It is obvious that $|G/Z(G)|=2^4=|G^{\prime}|^2$. \hfill$\Box$

\vspace{.2in} \noindent {\bf Acknowledgement.} Financial support of Council of Scientific and Industrial
Research, Government of India for the research of 
second and third author is gratefully acknowledged.

\end{document}